\documentclass[a4paper]{article}

\usepackage{amsmath,amsthm}
\usepackage{amsfonts,amssymb}
\usepackage{CJK}
\usepackage{graphicx}

\begin{CJK*}{GBK}{}

\end{CJK*}

\title{\bf  The Non-integrability of a Silnikov Equation
\footnote{ This work is supported by the fundamental
research funds for the central universities (12MS87).} }

\author{ Yanxia Hu
\thanks{Authors for correspondence. Email: yxiahu@163.com}
\date{\small \it School of Mathematics and Physics, North China Electric Power University,
Beijing,  102206 China } }

\begin{document}

\begin{CJK*}{GBK}{}
\CJKtilde \CJKfamily{song} \maketitle

\noindent{\small {\bf Abstract} \quad Based on the Lie group
theory, the one-parameter Lie group admitted by a Silnikov
equation (see \cite{Guan1}) is studied. The result reveals that
the Silnikov equation accepts no global analytical non-trivial
one-parameter Lie group. In this sense, the Silnikov equation is
not integrable in quadrature. \vskip 2mm {\bf Key words} \quad
Silnikov equation, Lie group, Infinitesimal generator
 \vskip 2mm

{\bf MR(1991)  Subject Classification} \quad 35B40; 35B45}


\vskip 5mm

\noindent {\large \bf 1 \quad Introduction}

\vskip 3mm

In the references \cite{Guan1, Guan2}, Keying Guan proposed and
studied the dynamical system described by
\begin{equation}\left\{\begin{array}{ccl}
  \displaystyle\dfrac{dx}{dt} &=& P(x,y,z)=y\\
  \displaystyle\dfrac{dx}{dt} &=& Q(x,y,z)=z\\
  \displaystyle\dfrac{dz}{dt} &=& R(x,y,z)=x^3-a^2x-y-bz,
\end{array}\right.\label{1}
\end{equation}
where $a,b$ are both positive constants, and he found system
(\ref{1}) is an important particular case of Silnikov equation.
Based on both qualitative method and numerical tests, he found
system (\ref{1}) has a series of interesting phenomena. For
example, this system may have "faint attractor (this term is
suggested by the author)", spatial limit closed orbits with
different rotation numbers and bifurcation of the limit closed
orbits. In \cite{Guan2}, Guan and Beiye Feng studied the
period-doubling cascades of the spatial limit closed orbit of
system (\ref{1}), and so on. It is revealed that the system is a
very rich source and an ideal model of three-dimensional
autonomous ordinary differential system in the research of the
bifurcation and chaos theory.

In this manuscript, we proceed to study the three-dimensional
system using Lie group theory. Professor Guan (the author of
\cite{Guan1}) has noticed that the system may not accept any
global analytical one-parameter Lie group excepting for a trivial
Lie group, and in this sense, one can conclude strictly that the
system is not integrable in quadrature. But he mentioned that he
has not given a detailed proof (see
http://blog.sciencenet.cn/u/guanky ).

In the fact, system (\ref{1}) is corresponding to a third-order
ordinary differential equation \begin{equation}a^2
y=y^3-y'''-y'-by'',\label{2}\end{equation} where $y'$ is the
derivative about variable $y$ with respect to variable $t.$ Since
the variable $t$ is not appeared explicit in equation (\ref{2}),
it admits a trivial one-parameter Lie group $t\rightarrow t+c.$
Let $z(y)=y',$ the above equation can turn to the following second
order equation
\begin{equation}z^2z''+zz'^2+z+a^2y-y^3+bz'z=0.\label{3}\end{equation}
The one-parameter Lie group admitted by (\ref{3}) is same as that
admitted by (\ref{2}) excepting for the trivial one-parameter Lie
group admitted by (\ref{2}). In this manuscript, we consider
system (\ref{3}) using Lie group theory, and find the fact that
the system accepts no global analytical one-parameter Lie group.

\vskip 5mm

\noindent {\large \bf 2 \quad  The One-parameter Lie group
admitted by the system}

\vskip 3mm

For the sake of convenience, we write equation (\ref{3}) as
\begin{equation}x^2\ddot{x}+x\dot{x}^2+x+a^2t-t^3+b\dot{x}x=0.\label{4}\end{equation}
We suppose that $V=\xi(t,x)\frac{\partial}{\partial
t}+\eta(t,x)\frac{\partial}{\partial x}$ is the infinitesimal
generator of Lie group admitted by equation (\ref{4}). Therefore,
$$V^{(2)}=\xi(t,x)\frac{\partial}{\partial
t}+\eta(t,x)\frac{\partial}{\partial x}
+\eta^{(1)}(t,x,\dot{x})\frac{\partial}{\partial \dot{x}}+
\eta^{(2)}(t,x,\dot{x},\ddot{x})\frac{\partial}{\partial
\ddot{x}}$$ be the 2th-extended infinitesimal generator, where
$$\eta^{(1)}(t,x,\dot{x})=\eta_t+(\eta_x-\xi_t)\dot{x}-\xi_x\dot{x}^2,$$
$$\eta^{(2)}(t,x,\dot{x},\ddot{x})=\eta_{tt}+(2\eta_{tx}-\xi_{tt})\dot{x}+
(\eta_{xx}-2\xi_{tx})\dot{x}^2-\xi_{xx}\dot{x}^3+(\eta_x-2\xi_t)\ddot{x}-3\xi_x\dot{x}\ddot{x}.$$
As we know (ref.{\cite{Blu}}), Equation (\ref{4}) will accept the
above symmetry group if and only if
\begin{equation}V^{(2)}(\ddot{x}-f(t,x,\dot{x}))=0 \; when  \; \ddot{x}=f(t,x,\dot{x}), \label{5}\end{equation}
where
$f(t,x,\dot{x})=\dfrac{t^3-a^2t-x-\dot{x}^2x-b\dot{x}x}{x^2}.$

Substituting
$\eta^{(1)}(t,x,\dot{x}),\eta^{(2)}(t,x,\dot{x},\ddot{x})$ and
$\ddot{x}=\dfrac{t^3-a^2t-x-\dot{x}^2x-b\dot{x}x}{x^2}$ to
(\ref{5}), we have
$$\begin{array}{lccc}\xi\dfrac{a^2-3t^2}{x^2}+\eta\dfrac{2(t^3-a^2t-x-\dot{x}^2x-b\dot{x}x)+(1+\dot{x}^2+b\dot{x})x}{x^3}\\
+[\eta_t+(\eta_x-\xi_t)\dot{x}-\xi_x\dot{x}^2]\dfrac{b+2\dot{x}}{x}\\
+[\eta_{tt}+(2\eta_{tx}-\xi_{tt})\dot{x}+(\eta_{xx}-2\xi_{tx})\dot{x}^2-\xi_{xx}\dot{x}^3\\
+(\eta_x-2\xi_t)\dfrac{t^3-a^2t-x-\dot{x}^2x-b\dot{x}x}{x^2}
-3\xi_x\dot{x}\dfrac{t^3-a^2t-x-\dot{x}^2x-b\dot{x}x}{x^2}]\\
=0\end{array}$$ After straightforward computing, it appears
$$\begin{array}{lccc}\dot{x}(-bx\eta+b(\eta_x-\xi_t)x^2-bx^2(\eta_x-2\xi_t)+2\eta_tx^2+(2\eta_{tx}-\xi_{tt})x^3-3\xi_x(t^3-a^2t-x)x)\\
+\dot{x}^2(2b\xi_xx^2-\eta x+2(\eta_x-\xi_t)x^2+(\eta_{xx}-2\xi_{tx})x^3-(\eta_x-2\xi_t)x^2)\\
+\dot{x}^3(-2\xi_xx^2-\xi_{xx}x^3+3\xi_xx^2)\\
+\xi(a^2-3t^2)x+2\eta(t^3-a^2t-x)+\eta
x+\eta_{tt}x^3+(\eta_x-2\xi_t)(t^3-a^2t-x)x+b\eta_tx^2\\
=0.\end{array}$$ The resulting determining equations for $\xi$ and
$\eta$ are given by
\begin{equation}\begin{array}{lccc}-bx\eta+b(\eta_x-\xi_t)x^2-bx^2(\eta_x-2\xi_t)+2\eta_tx^2+(2\eta_{tx}-\xi_{tt})x^3-3\xi_x(t^3-a^2t-x)x=0,\\
2b\xi_xx^2-\eta
x+2(\eta_x-\xi_t)x^2+(\eta_{xx}-2\xi_{tx})x^3-(\eta_x-2\xi_t)x^2=0,\\
-2\xi_xx^2-\xi_{xx}x^3+3\xi_xx^2=0,\\
\xi(a^2-3t^2)x+2\eta(t^3-a^2t-x)+\eta
x+\eta_{tt}x^3+(\eta_x-2\xi_t)(t^3-a^2t-x)x+b\eta_tx^2
=0.\end{array}\label{6}\end{equation} Form the above third
equation of (\ref{6}), one sees that
\begin{equation}\xi=\dfrac{1}{2}f_1(t)x^2+f_2(t),\label{7}\end{equation} where $f_1(t)$ and $f_2(t)$ are analytical functions. Then, we consider the
second equation of (\ref{6}). After simplifying the equation, we
are led to
$$-\eta+2b\xi_xx+\eta_xx+(\eta_{xx}-2\xi_{tt})x^2=0.$$
Substituting $\xi$ to the above equation, we have
$$\eta_{xx}x^2+\eta_xx-\eta=2f'_1(t)x^3-2bf_1(t)x^2.$$ Without lost of generality, we can let $\eta=g_1(t)x^3+g_2(t)x^2+g_3(t)x+g_4(t),$
where $g_i(t), i=1,2,3,4$ are all analytical functions. According
to
the above equation, one can obtain easily $$\begin{array}{lcc}g_1(t)=\dfrac{1}{4}f'_1(t),\\
g_2(t)=-\dfrac{2b}{3}f_1(t),\\
g_4(t)=0.
\end{array}$$
Thus far,
\begin{equation}\eta=\dfrac{1}{4}f'_1(t)x^3-\dfrac{2b}{3}f_1(t)x^2+g_3(t)x.\label{8}\end{equation}
Next, we consider the first equation of (\ref{6}). Substituting
(\ref{7}) and (\ref{8}) to the first equation of (\ref{6}), one
has
\begin{equation}\begin{array}{lcc}2f_1''(t)x^4
-\dfrac{15}{4}f_1'(t)x^3
+(\dfrac{2b^2+9}{3}f_1(t)+4g'_3(t)-f''_2(t))x^2\\
+(-3f_1(t)(t^3-a^2t)-bg_3(t)+bf'_2(t))x=0.\end{array}\label{9}\end{equation}
The left side of (\ref{9}) is a polynomial on $x.$ So, we can
obtain that $f_1(t)$ is a constant function $c_1$ and
\begin{equation}\begin{array}{lcc}\dfrac{2b^2+9}{3}f_1(t)+4g'_3(t)-f''_2(t)=0,\\
-3f_1(t)(t^3-a^2t)-bg_3(t)+bf'_2(t)=0.\end{array}\label{10}\end{equation}
After integrating the first equation of (\ref{10}) and multiplying
$b$, it appears
\begin{equation}4bg_3(t)-bf_2'(t)=-\dfrac{(2b^2+9)bc_1}{3}t+c_2b,\label{11}\end{equation}
where $c_2$ is an integrating constant. Adding the second equation
of (\ref{10}) to (\ref{11}), we have
$$g_3(t)=\dfrac{c_1(t^3-a^2t)}{b}-\dfrac{(2b^2+9)c_1t}{9}+\dfrac{c_2}{3}.$$
Based on (\ref{10}), we can obtain
$$f_2(t)=\dfrac{c_1t^4}{b}-\dfrac{2c_1a^2}{b}t^2-\dfrac{(2b^2+9)c_1t^2}{18}+\dfrac{c_2t}{3}+c_3,$$
where $c_3$ is an integrating constant. So,
$$\begin{array}{lcc}\xi=\dfrac{1}{2}c_1x^2+\dfrac{c_1t^4}{b}-\dfrac{2c_1a^2}{b}t^2-\dfrac{(2b^2+9)c_1t^2}{18}+\dfrac{c_2t}{3}+c_3,\\
\eta=-\dfrac{2b}{3}c_1x^2+(\dfrac{c_1(t^3-a^2t)}{b}-\dfrac{(2b^2+9)c_1t}{9}+\dfrac{c_2}{3})x.
\end{array}$$
Substituting $\xi$ and $\eta$ to the fourth equation of (\ref{3})
and  computing straightway, we can obtain $c_1=c_2=c_3=0.$ It
reveals $\xi(t,x)=0, \eta(t,x)=0.$ Consequently, equation
(\ref{3}) or (\ref{4}) admits a one-parameter Lie group with
infinitesimal generator $$V=\xi(t,x)\dfrac{\partial}{\partial
t}+\eta(t,x)\dfrac{\partial}{\partial x},$$ where $\xi(t,x)=0,
\eta(t,x)=0.$ That is, equation (\ref{3}) admits no global
analytical one-parameter Lie group. Equation (\ref{2}) admits no
global analytical one-parameter Lie group excepting for a trivial
Lie group $t\rightarrow t+c,$ and it is not integrable in
quadrature.

Clearly, when $b=0,$ (\ref{2}) is the special case
\begin{equation}a^2 y=y^3-y'''-y'.\label{12}\end{equation}The
above discussion can also indicate (\ref{12}) not admitting any
global analytical one-parameter Lie group excepting for a trivial
Lie group $t\rightarrow t+c,$ and it is not integrable in
quadrature.

On the general Silnikov equation mentioned in \cite{Guan1}, it can
be studied with the similar way. Relevant results will appear in
later manuscript.

\vskip 5mm

\noindent {\large \bf 3 \quad  Conclusions }

\vskip 3mm

In this manuscript, we have proved that the Silnikov equation is
not integrable in quadrature based on Lie group theory. So far, we
have seen the Silnikov equation possessing some dynamical
characters and particularity. It is believed that the Silnikov
equation still needs to be explored using other methods.

\vskip 3mm

\noindent {\large \bf \quad  Acknowledgments}

\vskip 3mm

The author would like to express her sincere gratitude to
professor Keying Guan (Beijing Jiaotong University) for his useful
and important suggestions on the method to proof the problem in
this note.

\end{CJK*}

\end{document}